\documentclass[reqno]{amsart}
\usepackage{amsmath}
\usepackage{amssymb}
\usepackage{amsthm}
\usepackage{graphicx,color,wrapfig}

\theoremstyle{plain} 
\newtheorem{lemma}{Lemma}

\newtheorem{prop}[lemma]{Proposition} 
\newtheorem{cor}[lemma]{Corollary}
\newtheorem{rem}[lemma]{Remark}

\newenvironment{Proof}[1][Proof]{
    \par
    \topsep 6pt plus 6pt
    \trivlist
    \item[\hskip\labelsep\bfseries #1.]\ignorespaces}%
    {\qed\endtrivlist
}

\newcommand{\R}{\mathbb{R}}

\begin{document}

\author{Hayk Mikayelyan}
\title[$\infty$-harmonic function, which is not $C^2$ on a dense 
subset]{An example of $\infty$-harmonic function, which is not $C^2$ on a dense subset}



\address{Hayk Mikayelyan \\
Mathematisches Institut \\
Universit\"at Leipzig\\\,\,\,
Augustusplatz 10/11\\
04109 Leipzig \\
Germany}
\email{hayk@math.uni-leipzig.de}

\maketitle

\begin{abstract}
We show that for certain boundary values McShane-Whitney's minimal-extension-like function 
is $\infty$-harmonic near the boundary and is not $C^2$ on a dense subset.
\end{abstract}

\vspace{5mm}

Let us consider the strip  $\{(u,v)\in\R^2:0<v<\delta\}$. 
The function we are going to construct
will be defined in this strip. Take a function $f\in C^{1,1}(\R)$ with 
$L_f:=\|f'\|_\infty$ and $L'_f:=Lip(f')$. 
Let us consider an analogue of the minimal extension of McShane and Whitney 
\begin{equation}
\label{McSh}
u(x,d):=\sup_{y\in\R}[f(y)-L|(x,d)-(y,0)|], 
\end{equation}
where $0<d<\delta$ and $L>L_f$. Note that in order to get the classical minimal 
extension of McShane and Whitney we have to take $L=L_f$.

From now on let us fix the function $f$, constants $L>L_f$ and $\delta>0$. 
We are going to find some conditions on $\delta>0$, which will make our 
statements to be true. The real number $x$ will be associated with 
the point $(x,\delta)\in\Gamma_\delta:=\{(u,v)\in\R^2:v=\delta\}$
and the real number $y$ with the point $(y,0)\in\Gamma_0$.
In the sequel the values of $u$ on the line $\Gamma_\delta$ will be 
of our interest and we write $u(x)$ for $u(x,\delta)$. 

\begin{prop}
\begin{equation}
\label{McSh1}
u(x)=\sup_{y\in\R}[f(y)-L\sqrt{\delta^2+(x-y)^2}]=
\max_{|y-x|\leq D\delta}[f(y)-L\sqrt{\delta^2+(x-y)^2}],
\end{equation}
where $D:=\frac{2LL_f}{L^2-L^2_f}$.
\end{prop}
\begin{Proof}
From the definition of $u$ we have that
$$
f(x)-L\delta\leq u(x)
$$
so it is enough to show that if $|x-y|>D\delta$ then
$$
f(y)-L\sqrt{\delta^2+(x-y)^2}< f(x)-L\delta.
$$
On the other hand from the boundedness of $f'$ we have
$$
f(y)-L\sqrt{\delta^2+(x-y)^2}\leq f(x)+L_f |x-y|-L\sqrt{\delta^2+(x-y)^2}.
$$
Thus we note that all values of $y$ for which 
$$
f(x)+L_f |x-y|-L\sqrt{\delta^2+(x-y)^2}<f(x)-L\delta
$$
can be ignored in taking supremum in the definition of $u$.

We write
$$
L_f |x-y|+L\delta<L\sqrt{\delta^2+(x-y)^2}
$$
and arrive at
$$
L_f^2|x-y|^2+2LL_f\delta|x-y|+L^2\delta^2<L^2\delta^2+L^2|x-y|^2
$$
Thus
$$
2LL_f\delta<(L^2-L_f^2)|x-y| \,\,\, <=>   \,\,\, |x-y|>D\delta.
$$
\end{Proof}

\vspace{4mm}

Let $y(x)$ be one of the points in $\{|y-x|\leq D\delta\}$, where the
maximum in (\ref{McSh1}) is achieved
\begin{equation}
\label{achieved}
u(x)=f(y(x))-L\sqrt{\delta^2+(x-y(x))^2}.
\end{equation}

\begin{figure}
\begin{center}
\input{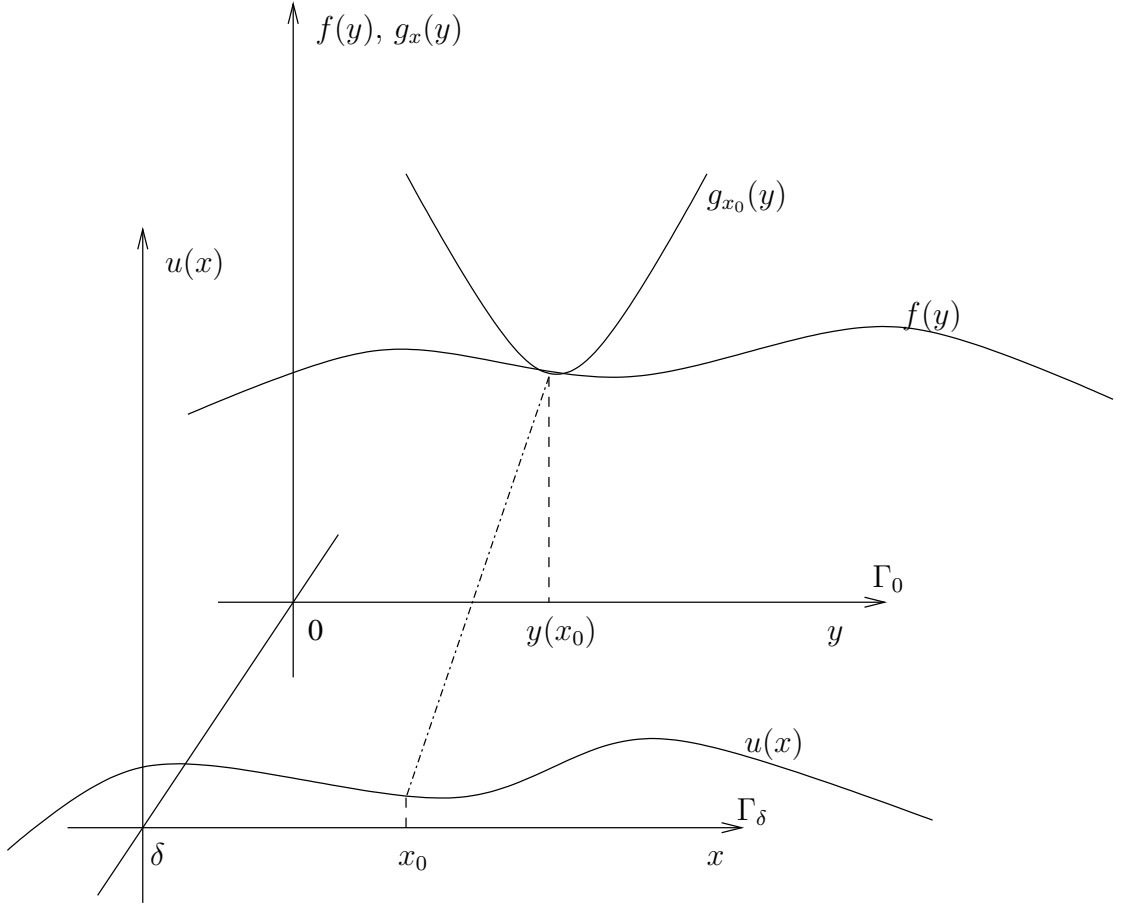}
\caption{Touched by hyperbola}
\label{cone}
\end{center}
\end{figure}

\begin{lemma}
If $\delta>0$ is small enough 
then for every $x\in\Gamma_\delta$ the point $y(x)$ is unique
and $y(x):\R\to\R$ is a bijective Lipschitz map.
\end{lemma}

\begin{Proof}
For each $x\in\Gamma_\delta$ consider the function 
$g_x(y):=u(x)+L\sqrt{\delta^2+(x-y)^2}$ defined on $\Gamma_0$.
The graph of $g_x$ is a hyperbola and the graph of any other function 
$g_{x'}$ can be obtained by a translation. 
Obviously $f(y)\leq g_x(y)$ on $\Gamma_0$ and $g_x(y(x))=f(y(x))$.
If at every point $y\in\Gamma_0$ the graph of $f$ can 
be touched from above by some hyperbola $g_x(y)$
then we will get the surjectivity of $y(x)$. 
To obtain this the following will be enough
\begin{equation}
\label{scdder}
g''_x(y)>L'_f,\,\,\,for \,\,\, all\,\,\,|y-x|\leq D\delta.
\end{equation}
So we arrive at 
\begin{equation}
\label{maincond}
\delta < \frac{L}{L'_f(1+D^2)^\frac{3}{2}}.
\end{equation} 
Note that also uniqueness of $y(x)$ follows from (\ref{scdder}); 
assume we have $y(x)$ and ${\tilde y}(x)$, then
$$
L'_f|y(x)-{\tilde y}(x)|<\left|\int_{y(x)}^{{\tilde y}(x)}g''_x(t)dt\right|=
|f'(y(x))-f'({\tilde y}(x))|\leq L'_f|y(x)-{\tilde y}(x)|.
$$
We have used here that
\begin{equation}
\label{firstder}
f'(y(x))=g_x'(y(x))=\frac{L(y(x)-x)}{\sqrt{\delta^2+(y(x)-x)^2}}
\end{equation}
(derivatives in $y$ at the point $y(x)$).  

The injectivity of the map $y(x)$ follows from differentiability of $f$.
Assume $y_0=y(x)=y({\tilde x})$, so we have 
$$
f(y_0)=g_x(y_0)=g_{{\tilde x}}(y_0).
$$
On the other hand $f(y)\leq min (g_x(y),g_{{\tilde x}}(y))$ this contradicts
differentiability of $f$ at $y_0$.

The monotonicity of $y(x)$ can be obtained using same arguments;
if $x<{\tilde x}$ then the `left' hyperbola $g_x(y)$ touches the 
graph of $f$ `lefter' than the `right' hyperbola $g_{{\tilde x}}(y)$, 
since both hyperbolas are above the graph of $f$.

Now we will prove that $y(x)$ is Lipschitz.
From (\ref{firstder}) it follows that
\begin{equation}
\label{MFEq}
y(x)-x=\frac{\delta f'(y(x))}{\sqrt{L^2-(f'(y(x)))^2}}.
\end{equation}
Taking $Y(x):=y(x)-x$ we can rewrite this as
\begin{equation}
\label{Banach}
Y(x)=\frac{\delta f'(Y(x)+x)}{\sqrt{L^2-(f'(Y(x)+x))^2}}=\delta\Phi(f'(Y(x)+x)),
\end{equation}
where $\Phi(t)=\frac{t}{\sqrt{L^2-t^2}}$.
For $ \delta < \frac{(L^2-L^2_f)^\frac{3}{2}}{L^2L'_f}$ we
can use Banach's fix point theorem and get that this functional
equation has unique continuous solution. On the other hand it is
not difficult to check that
$$
\left|\frac{Y(x_2)-Y(x_1)}{x_2-x_1}\right|\leq \frac{\delta C}{1-\delta C},
$$
where $C=\frac{L^2L'^2_f}{(L^2-L^2_f)^\frac{3}{2}}$.
\end{Proof}

\begin{cor}
If $\delta$ is as small as in the previous Lemma, then the function
$u$ is $\infty$-harmonic in the strip  between $\Gamma_0$ and $\Gamma_\delta$.
\end{cor}

\begin{Proof}
This follows from the fact that if we take the strip with boundary values
$f$ on $\Gamma_0$ and $u$ on $\Gamma_\delta$ then McShane-Whitney's minimal and maximal
solutions will coincide, obviously with $u$.
\end{Proof}

\begin{rem}
We can rewrite (\ref{MFEq}) in the form
\begin{equation}
\label{MFEq1}
x(y)=y-\frac{\delta f'(y)}{\sqrt{L^2-(f'(y))^2}},
\end{equation} 
where $x(y)$ is the inverse of $y(x)$. 
This together with (\ref{achieved}) gives us the following 
$$
u(x(y))=f(y)-\frac{\delta L^2}{\sqrt{L^2-(f'(y))^2}}.
$$
Using the recent result of O.Savin that $u$ is $C^1$, we conclude 
that function $x(y)$ is as regular as $f'$,
so we cannot expect to have better regularity than Lipschitz.
\end{rem}

\begin{lemma}
If $\delta>0$ is as small as above and function $f$ is not twice
differentiable at $y_0$, then the function $u$ is not twice
differentiable at $x_0:=x(y_0)$. 
\end{lemma}

\begin{Proof}
First note that for all $x$ and $y$, such that $x=x(y)$
we have 
$$
u'(x)=f'(y).
$$
This can be checked analytically but actually is a trivial
geometrical fact; the hyperbola 'slides' in the direction of the
growth of $f$ at point $y$, thus the cone which generates this hyperbola 
and 'draws' with its peak the graph of $u$ moves in same direction which is 
the direction of the growth of $u$ at point $x=x(y)$.

Now assume we have two sequences $y_k\to y_0$ and 
${\tilde y}_k\to y_0$ such that
$$
\frac{f'(y_k)-f'(y_0)}{y_k-y_0}\to\underline{f''}(y_0) \,\,\,and\,\,\,
\frac{f'({\tilde y}_k)-f'(y_0)}{{\tilde y}_k-y_0}\to\overline{f''}(y_0)
$$
and $\underline{f''}(y_0)<\overline{f''}(y_0)$.
Let us define appropriate sequences on $\Gamma_\delta$ denoting by 
$x_k:=x(y_k)$ and by ${\tilde x}_k:=x({\tilde y}_k)$ and 
compute the limits of
$$
\frac{u'(x_k)-u'(x_0)}{x_k-x_0} \,\,\,and\,\,\,
\frac{u'({\tilde x}_k)-u'(x_0)}{{\tilde x}_k-x_0}.
$$
We have 
$$
\frac{u'(x_k)-u'(x_0)}{x_k-x_0}=\frac{f'(y_k)-f'(y_0)}{y_k-y_0}\frac{y_k-y_0}{x_k-x_0}
$$
the first multiplier converges to $\underline{f''}(y_0)$, let us compute the limit of 
the second one. From (\ref{MFEq1}) we get that
$$
\frac{x_k-x_0}{y_k-y_0}\to 1-\delta \Phi'(f'(y_0))\underline{f''}(y_0),
$$
where $\Phi(t)=\frac{t}{\sqrt{L^2-t^2}}$.
Thus
$$
\frac{u'(x_k)-u'(x_0)}{x_k-x_0}\to
\frac{\underline{f''}(y_0)}{1-\delta \Phi'(f'(y_0))\underline{f''}(y_0)},
$$
and analogously
$$
\frac{u'({\tilde x}_k)-u'(x_0)}{{\tilde x}_k-x_0}\to
\frac{\overline{f''}(y_0)}{1-\delta \Phi'(f'(y_0))\overline{f''}(y_0)}.
$$
To complete the proof we need to use the monotonicity of the function
$$
\frac{t}{1-\delta C t}, \,\,\,-L'_f<t<L'_f,
$$
where $\frac{1}{L}<C<\frac{L^2}{(L^2-L^2_f)^\frac{3}{2}}$.
\end{Proof}

We would like to note that if the function $f$ is not $C^2$ at a point $y$ then
$u$ constructed here is not $C^2$ on the whole line connecting $y$ and $x(y)$.
So choosing $f$ to be not twice differentiable on a dense set we can get a function
$u$ which is not $C^2$ on the collection of corresponding line-segments.
Note that a similar example is the distance function from a convex 
set, whose boundary is $C^1$ and not $C^2$ on a dense subset.
Then the distance function is $\infty$-harmonic and is not $C^2$ on appropriate lines.

Our example has the property of having constant $|\nabla u|$ on gradient flow curves
(lines in our case). It would be interesting to find a general answer to the question: 

\vspace{3mm}

\noindent 'What geometry do the gradient flow curves have, on which $|\nabla u|$ is not constant?'

\vspace{3mm}

From Aronsson's results we know that $u$ is not $C^2$ on such a curve.
This is our motivation for the investigation of $C^2$-differentiability
of $\infty$-harmonic functions.

\vspace{5mm}

The author has only one item in the list of references. The 
history and the recent developments of the theory of $\infty$-harmonic 
functions, as well as a complete reference list could be found in that paper.

\vspace{5mm}

\subsection*{Acknowledgment}
The author is grateful to Gunnar Aronsson, Michael Crandall and 
Arshak Petrosyan for valuable discussions.

\end{document}